\documentclass[11pt,reqno, a4paper]{amsart}

\usepackage{amsmath,amsthm}
\usepackage{amssymb}
\usepackage{amsfonts}
\usepackage{amscd}
\usepackage{enumerate}
\usepackage{color}
\usepackage[all]{xy}
 at 9truept

\newtheorem{thm}{Theorem}[section]
\newtheorem{theorem}[thm]{Theorem}

\newtheorem{corollary}[thm]{Corollary}

\newtheorem{proposition}[thm]{Proposition}
\newtheorem{definition}[thm]{Definition}

\theoremstyle{definition}
\newtheorem{example}[thm]{Example}
\newtheorem{examples}[thm]{Examples}
\newtheorem{remark}[thm]{Remark}
\newtheorem{remarks}[thm]{Remarks}

\newtheorem{notation}[thm]{Notation}
\newtheorem{free text}[thm]{}

\numberwithin{equation}{section}

\newcommand{\N} {\mathbb{N}}
\newcommand{\Z} {\mathbb{Z}}

\newcommand{\F}{{\mathcal F}}
\newcommand{\e} {\epsilon}

\newcommand{\Der}{\operatorname{Der}}

\newcommand{ \strad} {{\rm strad}}

\newcommand{\h}{{\mathfrak h}}
\newcommand{\g}{{\mathfrak g}}
\newcommand{\Co}{{\rm Coind}}
\newcommand{\Ind}{{\rm Ind}}

\newcommand{{\Uop}}{{U^{\rm op}}}
\newcommand{{\Ae}}{{A^{\rm e}}}
\newcommand{{\Aop}}{{A^{\rm op}}}
\newcommand{{\Bop}}{{B^{\rm op}}}

\newcommand{{\op}}{{{\rm op}}}
\newcommand{{\coop}}{{{\rm coop}}}
\newcommand{{\sop}}{{*^{\rm op}}}
\newcommand{{\Ber}}{{{\rm Ber}}}

\newcommand{\lact}{{\,\raise1pt\hbox{$\scriptscriptstyle{\rhd}$}\,}}                  %
\newcommand{\ract}{{\,\raise1pt\hbox{$\scriptscriptstyle{\lhd}$}\,}}                  
\newcommand{\blact}{{\,\raise1pt\hbox{$\scriptscriptstyle{\blacktriangleright}$}\,}}  %
\newcommand{\bract}{{\,\raise1pt\hbox{$\scriptscriptstyle{\blacktriangleleft}$}\,}}   %

\begin{document}

\title{Duality properties for induced and coinduced representations in positive characteristic}
\author{Sophie Chemla\\
Sorbonne Université, Université  Paris Cité, CNRS, IMJ-PRG, F-75005 Paris, France.\\
{\it e-mail} : sophie.chemla@sorbonne-universite.fr}


\begin{abstract}
Let $k$ be a field of positive characteristic $p>2$. 
Generalizing a result of \cite{F-S}, we study the links between  coinduced representations and induced representations in the case of restricted Lie superalgebras. As a corollary, 
we prove a duality property concerning the kernel of coinduced representations of Lie $k$-superalgebras. This property was already proved by M. Duflo (\cite{Du})
for Lie algebras in any characteristic under more restrictive finiteness conditions.  It was then generalized to Lie superalgebras in characteristic 0 in  previous works (\cite{C0}, \cite{C1}, \cite{C4}). 
\end{abstract}

\maketitle

\section{Introduction}

Assume that  $\g=\g_{\overline{0}}\oplus \g_{\overline{1}}$ is a Lie superalgebra over a field $k$ 
of characteristic $p>2$ and $\h \subset \g$ is a Lie subsuperalgebra of $\g$. From a 
representation $(\pi, V)$ of $\h$, one can construct a representation of $\g$ in two ways: 
\begin{itemize}
\item Induction: $\Ind_\h^\g (V)=U(\g)\otimes_{U(\h)} V$ with left  $U(\g)$-module structure given by left multiplication;
\item  Coinduction: $\Co_\h^\g (\pi )=Hom_{U(\h)}(U(\g),V)$ with left $U(\g)$-module structure given by the transpose of right multiplication.
\end{itemize}

It is easy to see that the contragredient representation of $\Ind_\h^\g(V)$ is isomorphic to the coinduced representation from the contragredient representation $\pi^*$ of $\pi$ (\cite{Di}).

M. Duflo (\cite{Du}) proved that, in any characteristic, for a finite dimensional Lie algebra $\g$, the kernel $I_\pi\subset U(\g)$ of 
$\Co_\h^\g (\pi)$ satisfies  the duality property 
\begin{equation}\label{duality kernels}
\check{I}_\pi=I_{\pi^* \otimes k_{-trad_{\g/\h}}}
\end{equation}
where $\check{(-)}$ is the antipode of $U(\g)$,  $trad_{\g/\h}$ is the character $tr\circ ad_{\g/\h}$ of $\h$ and 
$k_{trad_{\g/\h}}$ the one dimensional representation of $\h$ it defines. 
In characteristic $0$, this duality property was extended to a Lie superalgebra $\g$ such that only $\g/\h$ is finite dimensional. In this case, the character $tr\circ ad_{\g/\h}$ is replaced by the character $str\circ ad_{\g/\h}$. 

The starting point of this article was to  treat the case of a Lie superalgebra with 
$\g/\h$ finite dimensional when $k$ is of  positive characteristic  $p>2$. We make use of a new ascending 
 filtration of $U(\g)$ (see \ref{definition filtration  mathcal F}) consisting in $U(\h)$-modules.


From now on, we assume that $\g/\h$ is finite dimensional and $k$ of positive characteristic $p>2$.  
In the case where $\g$ is a restricted Lie $k$-superalgebra with restricted enveloping superalgebra 
$U^\prime (\g)$ and $\h$ is a  restricted subsuperalgebra of $\g$,  stronger results holds. 
For a restricted 
representation $(\pi, V)$ of $\h$, we introduce the restricted induced representation from $\pi$ and the restricted 
coinduced representation from $\pi$: 
\begin{itemize}
\item Restricted induction: $\Ind_{U^\prime (\h)}^{U^\prime (\g)} (V)=U^\prime (\g)\otimes_{U^\prime (\h)} V$ 
with left  $U^\prime(\g)$-module structure given by left multiplication;
\item  Restricted coinduction: $\Co_{U^\prime(\h)}^{U^\prime(\g)} (\pi )=Hom_{U^\prime(\h)}(U^\prime(\g),V)$ with left 
$U^\prime (\g)$-module structure given by the transpose of right multiplication.
\end{itemize}

Generalizing a result of Borho-Brylinski (\cite{B-B}), it was proved in \cite{C4} that, in characteristic 0, the induced representation of a Lie superalgebra could be realized in terms of Grothendieck local cohomology. An analog result holds for restricted Lie superalgebras where local cohomology with coefficients in 
$\Co_{U^\prime(\h)}^{U^\prime(\g)} (\pi )$ is concentrated in degree $0$. It was proved in \cite{F-S} for finite dimensional Lie algebras:\\\\
{\bf Theorem} \ref{Ind and Coind} {\it Let $k$ be a field of characteristic $p>2$. Assume that the $k$-superspace 
 ${\g}/{\h}$ is finite dimensional. Set $m=dim \left (\g/\h\right )_{\overline{1}}$ and denote by 
 $\Pi$ the functor "change of parity". Let $(\pi, V)$ be a representation of $U^\prime (\h)$. Denote by 
 $\Pi ^m k_{strad_{\g/\h}}$ the one dimensional representation of $U^\prime (\h)$ with character $strad_{\g/\h} $. 
The restricted induced representation $\Ind_{U^\prime (\h)}^{U^\prime (\g)}(\pi \otimes \Pi^m k_{strad_{\g/\h}})$ is isomorphic to $\Co_{ U^\prime (\h)}^{U^\prime (\g)}(\pi)$.}\\

In \cite{C0}, it was noticed that (for any characteristic) Berezin integral provides  a $\g$-invariant duality between $\Co_{ \h}^\g(\pi)$ and  
$\Co_{ \h }^\g (\pi^* \otimes \Pi^m k_{-strad_{\g /\h}}) $ in case the Lie superspace $\g/\h$ is totally odd. We extend this result to any restricted Lie superalgebra $\g$ in the positive characteristic case: \\ 

{\bf Theorem} \ref{definition of  Psi} {\it Let $k$ be a field of characteristic $p>2$. Assume that the $k$-superspace 
 ${\g}/{\h}$ is finite dimensional with odd dimension $m$. There exists   a  non degenerate $\g$-invariant duality $\Phi$ between 
 $\Co_{U^\prime( \h) }^{U^\prime (\g)} (\pi) $ and 
 $\Co_{U^\prime( \h)}^{U^\prime(\g)}(\pi^* \otimes  \Pi^m k_{-strad_{\g /\h}})$.} \\
 
 Theorem \ref{Ind and Coind} and Theorem \ref{definition of  Psi} are linked by the $U^\prime (\g)$-module isomorphism  
 
$$\begin{array}{rcl}
\Theta : Coind_{U^\prime (\h)}^{U^\prime (\g)} (\pi^*)&\to &
 {\mathcal I}nd_{U^\prime (\h)}^{U^\prime (\g)} (\pi)^*\\
\lambda & \mapsto & \left [ u \otimes v\mapsto <\lambda(\check{u}),v>\right ].\\
\end{array}$$
 
We extend the duality property (\ref{duality kernels}) to Lie superalgebras as a corollary of  Theorem \ref{Ind and Coind}.

\section{Notation and preliminaries}

In this article, $k$ will be a commutative field of characteristic $p$. For most definitions about supermathematics, we refer the reader to \cite{Leites}. 
We will denote by $\overline{0}$ and $\overline{1}$ the elements of $\Z/2\Z$. We will call superspace a $k$-vector space graded over $\Z/2\Z$, 
 $V=V_{\overline{0}}\oplus V_{\overline{1}}$. If $v\in V_{i}$,  its degree will be denoted $\mid v \mid =i$. 
 As usual, formulas are meant for homogeneous elements and extended to any element by linearity. Let $V$ and $W$ be two superspaces. If $f$ is a morphism of degree $i$ from $V$ to $W$and if $v$ is in $V_j$, the element $f(v)$ will be also denoted $<f,v>$ and 
 we will use the notation
 $$<v,f>=(-1)^{ij}f(v)$$
especially  when it avoids the appearance of signs. If $V$ is a superspace, one defines the superspace $\Pi V$ which,  as a vector space, is equal to $V$ but the grading of which is 
 $(\Pi V)_{\overline{0}}=V_{\overline{1}}$ and $(\Pi V)_{\overline{1}}=V_{\overline{0}}$ . Let us introduce the map $\pi :V \to \Pi V$ which, as a morphism of vector spaces,  equals identity. It is of degree $\overline{1}$.  The functor $\Pi$ is called functor "change of parity". The symmetric superalgebra of $V$ will be denoted $S(V)$. 
 
 Let $A$ be an associative supercommutative superalgebra with unity and $M$ be an $A$-module. A basis of $M$ is a family 
 $(m_i)_{i \in I \amalg J} \in M_{\overline{0}}^I \times M_{\overline{1}}^J$ such that each element of $M$ can be expressed in a unique way 
 as a linear combination of $(m_i)_{i \in I \amalg J}$. If $I$ and $J$ are finite, their cardinality is independent of the basis of the $A$-module $M$. Then, the dimension of $M$ is 
 $(\mid I \mid , \mid J \mid )\in \N^2$. If $(e_1, \dots , e_n )$ is a basis of the $A$-module $M$, then the family $(e^1, \dots ,e^n)$ or $(e_1^*, \dots ,e_n^*)$ where 
 $<e_i, e^j>=\delta_{i,j}$ is a basis of $Hom_A (M,A)$ called the dual basis of $(e_1, \dots , e_n)$. Moreover, if $M$ is an $A$-module, then $\Pi M$ has a natural 
 $A$-module structure defined by:
 $$\forall m \in M, \; \forall a \in A, \quad a\cdot \Pi m=(-1)^{\mid a \mid} \Pi (a\cdot m).$$
 The following proposition is proved in \cite{Manin} p.172.
 \begin{proposition}\label{Berezinian}
 Assume that $M$ is a free $A$-module with finite dimension  $dim M=(n,m)$. Set $M^*=Hom_A(M,A)$. 
 If $(e_1, \dots , e_{n+m})\in M_{\overline{0}}^n \times M_{\overline{1}}^m$ is a basis of the $A$-module $M$,  denote by 
 $d$  left multiplication by 
 ${\displaystyle \sum_{i=1}^{n+m}(-1)^{\mid e_i \mid +1}\Pi e_i \otimes e^i}$ in the superalgebra 
 $S_A(\Pi M\oplus M^*)$. The endomorphism $d$ does not depend on the choice of a basis. 
 The complex 
 $$J(M)^\bullet=\left (S_A(\Pi M\oplus M^*)=\oplus_{n\in \N}S^n _A(\Pi M)\otimes_A S_A(M^*), d\right )$$
 has no cohomology except in degree $n$. The $A$-module $H^{n}(J(M))$ is free of dimension (1,0) or (0,1). More precisely the element 
 $\Pi e_1 \dots \Pi e_{n} \otimes e^{n+1} \dots e^{n+m}$ is a cycle the class of which is a basis of $H^n(J(M))$. 
  \end{proposition}
  The module $H^n(J(M))$ is called the Berezinian module of $M$ and is denoted $Ber (M)$. The Berezinian module 
  generalizes the maximal wedge, which does not exist if $M_{\overline{1}}\neq \{0\}$.\\
  
  Denote by ${\mathfrak g}l(M)$ the Lie superalgebra of endomorphims  of $M$. It acts on $S_A(\Pi M\oplus M^*)$ and its action commute with the differential $d$. Thus, it acts on 
  $Ber(M)$ by a character called "supertrace" and denoted $str$.\\

  We will make use of the Lie derivative:

 \begin{definition}\label{Lie derivative}
 (\cite{C2})  Let $A$ be a supercommutative $k$-superalgebra such that $Der(A)$ is a  finitely generated free  $A$-module.  Set 
 $Der(A)^*=Hom_A(Der(A),A)$.  The adjoint action of the Lie superalgebra $Der (A)$ 
on the complex $J(Der(A)^*)^\bullet$ induces an action of $Der(A)$ on  $Ber (Der(A)^*)$ 
(with the notation of Proposition \ref{Berezinian}). If $D \in Der (A)$, the Lie derivative $L_D$ of $D$ is defined by 
$$\forall \omega \in Ber (Der(A)^*), \quad L_D(\omega)=D \cdot \omega.$$
\end{definition}

\begin{example}
If $A=k[X_1, \dots ,X_{m+n}]$ is a polynomial superalgebra with 
even variables $X_1, \dots , X_n$ and odd variables $X_{n+1}, \dots , X_{n+m}$. Denote by 
$\omega \in H^n\left ( J(Der (A)^*)\right )$ the class of 
$\Pi \left ( \frac{\partial }{\partial X_1}\right )^*\dots \Pi \left ( \frac{\partial }{\partial X_n}\right )^*
\frac{\partial }{\partial X_{n+1}}\dots \frac{\partial }{\partial X_{m+n}}$ in $J\left ( Der (A)^*\right )$.
 If $D=\sum_{i=1}^{n+m}f_i\frac{\partial }{\partial X_i}$, we define the divergence of $D$ by 
 $Div_D=\sum_{i=1}^{n+m}(-1)^{\mid f_i\mid \mid X_i\mid }\dfrac{\partial f_i}{\partial X_i}$.
The following assertion is easy to check:
$$D\cdot \omega =-Div (D)\omega.$$
\end{example}
  
  If $\g$ is a $k$-Lie superalgebra, we will write $U(\g)$ for  its enveloping superalgebra and $\Delta$ for the coproduct in $U(\g)$. 
  The Hopf superalgebra $U(\g)$ is filtered by the standard filtration $\left (F_nU(\g)\right )_{n \in \N}$.
  
  \begin{equation}\label{standard filtration}
 \left \{ \begin{array}{l}
  F_0U(\g)=k,\\
  \forall n\in \N^*,\quad F_nU(\g)=F_{n-1}U(\g)\cup Vect ( X_{1} \dots X_{n}, \quad  \forall i\in [1,n],\; X_i\in \g).
  \end{array}\right.
  \end{equation}
  If $V$ is a left $U(\g)$-module, then $V^*$ will be the contragredient module. Let us now describe the primitive elements of the Hopf superalgebra $U(\g)$. The following result is well known but, as  we did not find any reference, we give a proof of it. Denote by $\Delta$ the coproduct of $U(\g)$.
  
  \begin{proposition}\label{primitive elements}
  Let $\g$ be a Lie  superalgebra. 
 Let $(e_1, \dots , e_n)$ be a basis of $\g_{\bar{0}}$ and $( \epsilon_1, \dots  , \epsilon_m)$ be a basis of $\g_{\bar{1}}$.  The vector superspace of primitive elements of  $U(\g)$ is generated by  $\{e_i^{p^j}, \epsilon_s,  \quad (i,j,s) \in [1,n] \times \N \times \{1,m\} \}$.
 \end{proposition}
 
 \begin{notation} Let us introduce the following notation $\N^{n,m}:= \N^n \times \{0,1\}^m$.  If $(\underline{\bf a}, \underline{\alpha} )\in \N^n \times \{0,1\}^m$, we set 
 $$e^{\underline{\bf a}}\epsilon^{\underline{\alpha}}=e_1^{a_1}\dots e_n^{a_n}\epsilon_1^{\alpha_1}\dots \epsilon_n^{\alpha_n} .$$
 \end{notation}

 {\it Proof of Proposition \ref{primitive elements}}:
 
 Clearly the elements $e_i^{p^j}$ and $\epsilon_s$ are primitive elements. The family 
 $\left ( e^{\underline{\bf a}}\epsilon^{\underline{\alpha}}\right )_{(\underline{\bf a}, \underline{\alpha})\in \N^{n,m}}$ is a basis of $ U({\mathfrak g})$ and 
 $\left ( e^{\underline{\bf a^\prime}}\epsilon^{\underline{\alpha}^\prime}
 \otimes e^{\underline{\bf a^{\prime \prime}}}\epsilon^{\underline{\alpha^{\prime \prime}}}\right ) 
 _{(\underline{\bf a^\prime},\underline{\alpha^\prime}), ( \underline{\bf a^{\prime \prime}}, \underline{\alpha^{\prime \prime}})\in \N^{n,m}}$ is a basis of
  $U({\mathfrak g})\otimes U({\mathfrak g})$. 
  
  Let $x=\sum x_{\underline{\bf a}, \underline{\alpha}} e^{\underline{\bf a}}\epsilon^{\underline{\alpha}}$ be a primitive element of $U(\g)$. 
  One has 
  $$\Delta(x)-x\otimes 1-1 \otimes x = \sum_{
  \begin{array}{l}
 (\underline{\bf a^\prime}, \underline{\alpha^\prime}) , (\underline{\bf a^{\prime \prime}}, \underline{\alpha^{\prime \prime}})\neq  (\underline{0}, \underline{0}),\\
  (\underline{\bf a}^\prime, \underline{\alpha^{\prime}}) +(\underline{\bf a}^{\prime \prime}, \underline{\alpha^{\prime \prime}})=(\underline{\bf a}, \underline{\alpha}),\\
  
  \end{array}}
  x_{\underline{\bf a}, \underline{\alpha}}
  e^{\underline{\bf a^\prime}}\epsilon^{\underline{\alpha^\prime}}
  \otimes e^{\underline{\bf a^{\prime \prime}}}\epsilon^{\underline{\alpha^{\prime\prime}}}=0.$$
  If the term $e^{\underline{\bf a}}\epsilon^{\underline{\alpha}}$ involves more than one $e_i$ or $\epsilon_i$,  then  
  $\Delta(e^{\underline{\bf a}}\epsilon^{\alpha})-e^{\underline{\bf a}}\epsilon^{\alpha}\otimes 1-1 \otimes e^{\bf \underline{a}}\epsilon^{\alpha}\neq 0$.

 Thus 
  $x$ can be written $x=\sum_{i\in [1,n]}x_{a_i} e_i^{a_i}+ \sum_{s\in [1,m]}x_{s} \epsilon_s$ with $x_{a_i}$ and $x_s$ in $k$. 
  Let us now show that all the $a_i$'s are a power of $p$. Assume, it is not the case for $a_i$. 
  Let $p^{t_i} \in \N$ such that $p^{t_i}<a_i<p^{t_i+1}$. Set $a_i=p^t+b_i$ with $b_i\in [1,p-1]$.
  One has 
  $$\begin{array}{rcl}
  \Delta(e^{a_i})&= &\Delta(e^{p^t})\Delta(e^{b_i})\\
  &=& (e^{p^t}\otimes1+1 \otimes e^{p^t})(e^{b_i}\otimes 1+b_ie^{b_i-1}\otimes e_i+\cdots )
  \end{array}$$
  and the term $b_ie^{b_i-1}\otimes e_i$ is non zero so that 
  $\Delta(x)-x\otimes 1-1 \otimes x =0$ is not zero. $\Box$\\
  
  Unadorned tensor products are tensor products over $k$. When not specified, the duals are taken for $k$-vector spaces.

\section{Generalities on Lie-Rinehart superalgebras}

\begin{definition}
Let $A$ be a  supercommutative $k$-superalgebra. A $k-A$-Lie-Rinehart (\cite{Rinehart})
superalgebra (with anchor $\rho$) is a triple 
$(L, [\;,\; ], \rho )$ such that 
\begin{enumerate}
\item $(L, [\; , \; ])$ is a $k$-Lie superalgebra;
\item $L$ is an $A$-module;
\item The anchor $\rho: L \to \Der (A)$ is a $k$-Lie  superalgebra morphism and an $A$-module morphism such that: For all $(D,\Delta ) \in L$ and all $a\in A$,
$$[D,a\Delta ]=\rho(D)(a)\Delta +(-1)^{\mid a \mid \mid D \mid} a [D, \Delta ].$$
\end{enumerate}
\end{definition}
\begin{examples}
\begin{enumerate}
\item If $A=k$, a Lie-Rinehart superalgebra is a $k$-Lie superalgebra.
\item The $A$-module $\Der (A)$ is a $k$-$A$-Lie-Rinehart superalgebra with anchor equal to $id$.
\item Assume that $\g$ is a $k$-Lie superalgebra given with a Lie superalgebra morphism $\sigma :\g \to \Der (A)$. Then the $A$-module $A \otimes \g$ is  endowed with a unique $k-A$ Lie- Rinehart superalgebra such that 
\begin{itemize}
\item The anchor $\rho : A \otimes \g$ is defined by: For all $a \in A$ and $X \in \g$, $\rho (a X)=a\sigma (X)$.
\item The Lie bracket on $A \otimes \g$ extends that of $\g$.
\end{itemize}
The Lie-Rinehart superalgebra constructed that way is called the crossed product of $A$ with $\g$ and is denoted 
$(A\sharp \g , \sigma )$ or just $A\sharp \g$ when there is no ambiguity.
\item Poisson superalgebras also gives rise to Lie-Rinehart superalgebras but we won't use them in this article 
(\cite{Hu})
\end{enumerate}
\end{examples}

Rinehart (\cite{Rinehart}) has associated an enveloping algebra to a Lie-Rinehart algebra. This notion generalizes the enveloping algebra of a Lie algebra.

\begin{definition}  \label{def_U_A(L)}
 Let  $ (L,\rho) $  be a $k$-$A$-Lie-Rinehart superalgebra.  The  {\sl universal enveloping superalgebra}  of  $ L $  is the  $ k $--superalgebra
  $$  U_A(L)  \; := \;  T_k^+(A \oplus L) \Big/ I  $$
where  $ T_k(A \oplus L) $  is the tensor  $ k $-superalgebra  over  $ \, A \oplus L \, $  and  $ I $  is the two sided ideal in  $ T_k^+(A \oplus L) $  generated by the elements
  $$  a \otimes b - a \, b, \;\;   \qquad  a \otimes \xi - a \, \xi , \;\;
 \qquad  \xi \otimes \eta -(-1)^{\mid \eta \mid \mid \xi \mid}\eta \otimes \xi - [\xi,\eta], \;\; 
 \qquad  \xi \otimes a -(-1)^{\mid a \mid \mid \xi \mid}a \otimes \xi -\rho (\xi)(a)  $$
for all  $ \, a, b \in A \, $,  $ \, \xi , \eta \in L \, $.
\end{definition}

\begin{remark}
The anchor endows $A$ with a left $U_A(L)$-module structure. 
\end{remark}

In this article, we will mostly be in the case of the  crossed product superalgebra $A \sharp U(\g)$ given by a   coinduced superalgebra.\\

 Let $\h$ be a subLie superalgebra of $\g$ and let $(\pi, V)$ be a representation of ${\mathfrak h}$. We set 
 $$ \Co_\h^\g (\pi )=\{ \mu : U(\g) \to V, \; \forall u\in U(\g), \; \forall H \in \h, \; 
 <Hu, \mu >=\pi (H) <u,\mu> \}.$$
 The coinduced representation from $\pi$ is a representation of $\g$ over the space 
 $\Co_\h^\g (\pi)$ defined by 
 $$\forall (u,v) \in U(\g)^2, \forall \mu \in \Co_\h^\g (\pi ), \quad 
 <v, u\cdot \mu ,>=<vu, \mu >.$$
 The action of $u\in U(\g)$  on $\Co_\h^\g (\pi)$ will be denoted $\delta_u^\pi$. 
 
 If $(\pi , V)$ is the trivial representation (just denoted $k$), the coproduct of $U(\g)$ allows to endow $\Co_\h^\g (k)$ with a 
 $k$-superalgebra structure: If $(\lambda , \mu) \in \Co_{\h}^\g (k)^2$,
 $$\forall u \in U(\g), \quad <u, \lambda \mu >=
\sum  (-1)^{\mid \lambda \mid \mid u_{(2)}\mid }<u_{(1)}, \lambda><u_{(2)}, \mu>\; {\rm with}\;
 \Delta (u)=\sum u_{(1)}\otimes u_{(2)}.$$
 The superalgebra $A:=\Co_\h^\g(k)$ is local with maximal ideal
 $${\mathfrak a}=\{\lambda \in A,\quad <\lambda ,1>=0\}.$$

 The action of $X \in \g$ on $A:=\Co_\h^\g (k )$ is given by a derivation denoted 
 $\delta_X  \in \Der \left (\Co_\h^\g (k) \right )$ (instead of $\delta_X^k$).

 We can thus perform the crossed product construction $(\Co_\h^\g (k)\sharp \g , \delta)$. 
  From now on, we will write $A \sharp U(\g)$ for $(A \sharp U(\g) , \delta )$.\\
 The coproduct on $U(\g)$ allows to endow $\Co_\h^\g (\pi)$ with a left $\Co_\h^\g (k)$-module structure and (adding the coinduced representation of $\g$) 
 $\Co_\h^\g (\pi )$ becomes a left $U_A(\Co_\h^\g (k) \sharp \g)$-module.

If $(\pi , V)$ is a representation of $\h$, then the induced representation from $(\pi, V)$ is the $U(\g)$-module structure given by left multiplication on the superspace 
$$Ind_\h^\g(V)=U(\g)\otimes_{U(\h)}V.$$

\section{Algebraic structure on the $\Co_\h^\g (k)$}

In this section, we study $\Co_\h^\g (\pi)$. The restricted coinduced case will follow easily.

\subsection{The symmetric superalgebra} 

Let $V=V_{\overline{0}}\oplus V_{\overline{1}}$ be a $k$-vector superspace of finite dimension.  
Let  $(e_1, \dots, e_n)$ be a basis of $V_{\overline{0}}$ and $(\epsilon_{1}, \dots , \epsilon_m )$ be a basis of $V_{\overline{1}}$.
If $(\underline{\bf a}, \underline{\alpha})=(a_1, \dots , a_n, \alpha_1, \dots , \alpha_m) \in \N^n\times \{0,1\}^m$, we set 
$$e^{\underline{\bf a}}\epsilon^{\underline{\alpha}}=e_{1}^{a_1}\dots e_{n}^{a_n}\e_{1}^{\alpha_1}\dots \e_{m}^{\alpha_m}.$$
 The monomials  $e^{\underline{\bf a}}\epsilon^{\underline{\alpha}}$ form a basis of the $k$-superalgebra $S(V)$. \\

 
 As $S(V)$ is a cocommutative Hopf superalgebra, its dual $S(V)^*$ is a supercommutative $k$-superalgebra. 
  Moreover, $S(V)^*$ is a local superalgebra  with maximal ideal 
 $${\mathfrak m}=\{ f\in S(V)^*, \quad <f,1>=0\}.$$

  Define the elements 
 $\mu_{e^{\underline{\bf a}}}$  and $\nu_{\epsilon^{\underline{\alpha}}}$    of $S(V)^*$ by 
 $$
\begin{array}{l} 
< e^{\underline{\bf b}}\epsilon^{\underline{\beta}} , \mu_{e^{\underline{\bf a}}}>=
 \delta _{\underline{\bf a}, \underline{\bf b}} \delta_{\underline{0}, \beta};\\
 < e^{\underline{\bf b}}\epsilon^{\underline{\beta}} , \nu_{\epsilon^{\underline{\alpha}}}>=
 \delta _{\underline{0}, \underline{\bf b}} \delta_{\underline{\alpha}, \underline{\beta}}.\\
 \end{array}$$
 
$\left (\mu_{e^{\underline{\bf a}}} \nu_{\epsilon^{\underline{\alpha}}} \right )$  is a basis of ${S(V)^*}_f$ (the elements of $S(V)^*$ with finite rank) and satisfies the following properties:
\begin{proposition} If $a_i, b_i \in \N$ and $\alpha_s, \alpha_t \in \{0,1\}$, then 
 $$\begin{array}{l}
\mu_{e_i^{a_i}} \times \mu_{e_i^{b_i}}=
\begin{pmatrix}a_i+b_i\cr a_i \end{pmatrix} \mu_{e_i^{a_i+b_i}};\\
(\nu_{\e_s})^2=0;\\
 \mu_{e_i^{a_i}} \times \mu _{e_j^{a_j}}=\mu _{e_i^{a_i}e_j^{a_j}}\; {\rm if}\; i<j;\\
   \nu_{\e_s^{\alpha_s}} \times \nu _{\e_t^{\alpha_t}}=
   (-1)^{\mid \alpha_s \mid \mid \alpha_t \mid }\nu _{\e_s^{\alpha_s}\e_t^{\alpha_t}}\; {\rm if}\; s<t.\\\\
 \end{array}$$
 \end{proposition}
 
 {\bf Notation:}
We will rather make use of the 
  elements $\mu_{e_i^{p^j}} $ and $\nu_{\epsilon_s}$ of $S(V)^*$. They will be respectively  denoted $\mu_{i,j}$ and  $\nu_s$.\\
 
 They have nice properties:
   \begin{proposition}
   \begin{enumerate}
  \item $(\mu_{i,j})^p=0$ and $(\nu_s)^2=0$.
  \item $(S(V)^*)_f=k\left [\mu_{i,j}, \; \nu_s ,\quad  i\in [1,n], j\in \N, s \in [1,m] \right ]$.
 \item  ${\mathfrak m}^p=\{0\}$. 
  \end{enumerate}
 \end{proposition}
 
 {\it Proof :}
 
 We will only prove the first assertion. The second one follows directly.
 
 Let us prove that $(\mu_{e_{i}^{p^j}})^p=0$. It follows from the equality (see Proposition \ref{primitive elements}) 
  $$<(\mu_{e_{i}^{p^j}})^p, e_i^{p^{j+1}}>=<(\mu_{e_i^{p^j}})^{\otimes p}, \Delta^p(e_i^{p^{j+1}})>=0.$$

As $\nu_s$ is an odd element, we have $\nu_s^2=0$. $\Box$\\


 \subsection{The superalgebra $\Co_{\mathfrak h}^{\mathfrak g}(k)$}

 Let $\g$ be a Lie $k$-superalgebra and $\h$ be a Lie subsuperalgebra of $\g$. 
 Given a supplement  ${\mathfrak p}$ of $\mathfrak h$ in $\mathfrak g$ and a basis 
 $\underline{e}=(e_1, \dots , e_n; \epsilon_1, \dots, \epsilon_m)$ of it, 
 we will use the following notation in $U(\g)$ that will turn out to be more convenient to study the coinduced representations:
  
 \begin{notation} \label{usual notation}
 $$\begin{array}{l}
 e_{i,j}={e_i}^{p^j}, \\ 
 \left [0,p-1\right ]^{k,m}=\left [0,p-1 \right ]^{k} \times \{0,1\}^m,\\
 {\rm  For}\; r \in \N,  \; (\underline {a}, \underline{\alpha}) \in [0,p-1]^{n(r+1),m}, \;
 e^{\underline{a}}\epsilon^{\underline{\alpha}}=
\left [ \prod_{  i \in [1,n]} \prod_{j \in [0,r]}e_{i,j}^{a_{i,j}} \right ]\prod_{s=1}^m \epsilon_s^{\alpha_s}.
 \end{array}.$$ 
 From now on, we set 
 $$A:=\Co_\h^\g(k).$$
 
 To study $A$, we will  make use of a filtration $\F$ on $U(\g)$  different from the usual one $F$.
 It   is defined by: 
 \begin{equation}\label{definition filtration  mathcal F}
\left \{ \begin{array}{l}
  \F_{-1}U(\g)=U(\h),\\
 {\rm If\;} r \in \N, \; {\F}_rU(\g)= 
  U(\h)Vect < 1; \;e^{\underline{a}}\epsilon^{\underline{\alpha}},  
  \quad (\underline{a}, \underline{\alpha})\in [0,p-1]^{n(r+1),m}>.
  \end{array}\right .
  \end{equation}
   \end{notation}
   \begin{remarks}\label{properties of mathcal F}
   \begin{enumerate}
   \item ${\F}_rU(\g)$ is a left $U(\h)$-module.
   \item ${\F}_rU(\g)$ is a filtered coalgebra: $\Delta {\mathcal F }_rU(\g) \subset {\displaystyle \oplus_{t+t^\prime =r} {\mathcal F }_tU(\g) \otimes 
  { \mathcal F }_{t^\prime}U(\g)}$.
   \end{enumerate}
   \end{remarks}

  We define a map
 $J_0 :  \Co_{\mathfrak h}^{\mathfrak g}(k) \to Hom_k \left (S({\mathfrak p}), k \right )$ by: For all $\lambda \in \Co_{\mathfrak h}^{\mathfrak g}(k)$, all $r \in \N$, all 
 $(\underline{a}, \underline{\alpha}) \in [0,p-1]^{n(r+1),m}$,
  
 $$<e^{\underline{a}}\epsilon^{\underline{\alpha}}, J_0(\lambda)>= 
 <e^{\underline{a}}\epsilon^{\underline{\alpha}}, \lambda>.$$
The following proposition is proved in \cite{C1}.
 \begin{proposition} 
 The map $J_0$ is an isomorphism of superalgebras.
 
   \end{proposition}
   
   If $(\pi, V)$ is a $U(\h)$-module, 
   we also define  an isomorphism 
 $J_\pi : \Co_{\mathfrak h}^{\mathfrak g}(\pi) \to Hom \left (S({\mathfrak p} \right ), V)$ by: For all $r \in \N$ and all 
 $(\underline{a}, \underline{\alpha}) \in [0,p-1]^{n(r+1),m}$ and $\lambda_\pi \in \Co_{\mathfrak h}^{\mathfrak g}(\pi)$, 
 $$\left \{ \begin{array}{l}
 <e^{\underline{a}}\epsilon^{\underline{\alpha}}, J_\pi (\lambda_\pi)>= 
 <e^{\underline{a}}\epsilon^{\underline{\alpha}},\lambda_\pi >;\\
 <1, J_\pi (\lambda_\pi)>= <1,\lambda_\pi >.
\end{array} \right .$$
 One has for any $f \in A$ and $\lambda_\pi \in \Co_{\mathfrak h}^{\mathfrak g}(\pi)$,
 $$J_\pi (f\lambda_\pi)=J_0(f)J_\pi (\lambda_\pi).$$

  The filtration ${\F}$ (\ref{definition filtration mathcal F}) we have introduced on $U(\g)$ induces a filtration  on the $k$-algebra $A$  as follows :  
 $$ \dots \F_r A \subset \F_{r-1}A \subset \dots \subset \F_1 A  \subset \F_0 A \subset A$$
  where 
  $$\F_rA=\{\lambda \in A, \; <\lambda, u>=0 \; {\rm if} \, u \in {\F}_{r-1}U(\g)\}.$$
  As ${\displaystyle \cap_{r \in \N} \F_rA=\{0\}}$, this filtration defines a Hausdorff topology on $A$. Moreover, $A$ is complete for this topology. \\
  
  Let us denote 
  ${\mathfrak a}={J_0}^{-1}({\mathfrak m})=\{\lambda \in A, \; <1,\lambda>=0\}$. 
  
  \begin{remarks}
  Denote by $F_kU(\g)$ the usual ascending filtration of $U(\g)$ and $F_kA$ the decreasing filtration it induces on $A$.
  
  \begin{enumerate}
  
  \item  ${\mathfrak m}^q$ is included in $F_qA$ but is not equal if $q\geq p$. Indeed, 
  ${\mathfrak m}^p$ equals $k[\nu_1, \dots , \nu_m]$ but not $F_pA$.
  
  \item $J_0(F_q A)=Vect (\mu_{\underline{e}^{\underline{\bf a}}}, \; \mid {\bf a} \mid \geq q)$.
  \end{enumerate}
  \end{remarks}
 
  \begin{notation} Associated to the choice of a basis of a supplement ${\mathfrak p}$ of 
  $\mathfrak h$ in $\mathfrak g$, the elements 
  $\mu_{i,j} \in S({\mathfrak p})^*$ and 
  $\nu_s \in  S({\mathfrak p})^*$ were defined earlier for $i\in [1,n]$, $j \in \N$, $s \in [1,m]$.

  The element 
  ${J_{0}}^{-1}\ \left ( \mu_{i,j}\right ) \in A$ will be  denoted 
  $ \eta_{e_{i}^{p^j}}$ or $\eta_{i,j}$. The element ${J_{0}}^{-1}\ \left ( \nu_{s}\right ) \in A$ will be  denoted 
  $ \zeta_s$.
  
 The derivation ${J_0}^{-1} \circ \dfrac{\partial}{\partial \mu_{i,j}} \circ J_0$ 
  (respectively ${J_0}^{-1} \circ \dfrac{\partial}{\partial \nu_{s}} \circ J_0$)
   of $A$ will  be denoted $\partial_{i,j}$ (respectively $\overline{\partial_s}$). 
    \end{notation}

The study of differential operators of $Hom \left (S({\mathfrak p}), V \right )$ is transferred to 
 $\Co_{\mathfrak h}^{\mathfrak g} (\pi)$ by $J_\pi$. 
The derivation 
 ${J_\pi }^{-1} \circ \dfrac{\partial}{\partial \mu_{i,j}}\circ J_\pi$ 
 (respectively ${J_\pi }^{-1} \circ \dfrac{\partial }{\partial \nu_{s}}\circ J_\pi$)
 will be denoted $\partial^\pi_{i,j}$ (respectively $\overline{\partial_s^\pi}$) or more simply 
 $\partial_{i,j}$ (respectively $\overline{\partial_s}$) when there is no ambiguity.\\

The $k$-superalgebra $A$ is the projective limit of $\dfrac{A}{{\mathcal F}_r A}$ defined by the transpose of the natural injection 
${\mathcal F}_rU(\g) \hookrightarrow {\mathcal F}_{r^\prime}U(\g)$ if $r\leq r^\prime$.
Let us introduce the following notation: 
$$A^{\leq r}:=k[\{\eta_{i,j}, \zeta_s, \quad i \in [1,n],  \; j \in [0, r], s\in [1,m]\}]. $$
Let $\lambda \in A$.  The class of $\lambda$ in $\dfrac{A}{{\mathcal F}_{r+1}A}$ has  a unique representant, denoted by $\lambda^{\leq r}$,
 that is a polynomial in the $(\eta_{i,j})_{ j \leq r}$ and $(\zeta)_{s \in [1,m]}$. The map
 $$\begin{array}{rcl}
 A^{\leq r} & \to & \dfrac{A}{{\mathcal F}_{r+1} A}\\
 \lambda & \mapsto & \lambda+ \F_{r+1} A.
 \end{array}$$
 is a superalgebra isomorphism and $A=\varprojlim \dfrac{A}{{\mathcal F}_r A}$.

 For $X\in \g$, we define the derivation $\delta_X^{\leq r}$ of $A^{\leq r}$ as follows :
 $$\forall \lambda \in A^{\leq r}, \quad  \delta_X^{\leq r}(\lambda):= \delta_X(\lambda)^{\leq r}.$$
 \begin{proposition}
 If $\lambda \in A$ and $\lambda^{\leq r}$ is the representant of the class of $\lambda$ in  $\dfrac{A}{{\mathcal F}_{r+1}A}$ that is a polynomial in the 
$(\eta_{i,j})_{ j \leq r}$'s and $\zeta_s$'s, one has 
$$\delta _X (\lambda)= \lim\limits_{\mathcal F}\delta_{X}^{\leq r} (\lambda^{\leq r})=
 \lim\limits_{\mathcal F} \left [ \delta_X (\lambda^{\leq r})\right ]^{\leq r}.$$
 $\delta_X$ is continuous for the topology defined by the filtration ${\mathcal F}$ as 
 $\delta_X({\mathcal F}_r)\subset {\mathcal F}_{r-1}$.
\end{proposition}
We write 
$$\delta_{X}^{ \leq r}=\sum_{i=1}^n \sum_{j=0}^{r}f_{i,j}(X)^{\leq r}\partial_{i,j}+
\sum_{s=1}^m g_s(X)^{\leq r}\overline{\partial}_s.$$

More generally, for $r \in \N $, we define 
$$\Co_\h^\g (\pi)^{\leq r}:=\{\lambda_\pi \in \Co_\h^\g (\pi), \; 
<e^{\underline{a}}\epsilon^{\underline{\alpha}},\lambda_\pi>=0\; { \rm if } \; 
e^{\underline{a}}\epsilon^{\underline{\alpha}} \notin {\mathcal F}_rU(\g)\}=
\dfrac{\Co_\h^\g(\pi)}{{\mathcal F}_{r+1}A\, \Co_\h^\g(\pi)}.$$
$\Co_\h^\g (\pi)$ is the projective limit of the $\Co_\h^\g (\pi)^{\leq r}$. 
If $v \in V$, we define the element $\widehat{v}_\pi \in \Co_\h^\g (\pi)$ such that 
\begin{equation}\label{definition widehat{v_pi}}
\left \{\begin{array}{l}
<1,\widehat{v}_\pi>=v.\\
<e^{\underline{a}}\epsilon^{\underline{\alpha}}, \widehat{v}_\pi>=0, 
\; \forall r \in \N, \; \forall (\underline{a}, \underline{\alpha})\in \N^{n(r+1)}\times \{0,1\}^m.
\end{array}\right .
\end{equation}
The element $X \in \g$ defines a differential operator of degree one, $\delta_X^{\pi, \leq r}$,  of $\Co^\g_{ \h}{(\pi)}^{\leq r}$. We will write 
$$\delta_X^{\pi, \leq r}=F_X^{\pi, \leq r} + \sum_{i=1}^n \sum_{j=0}^r f_{i,j}^{\leq r}(X)\partial_{i,j}
+\sum_{s=1}^m g_s (X)^{\leq r}\overline{\partial_s}$$
where $F_X^{\pi, \leq r}$ is the element of $End_{A^{\leq r}}[\Co_\h^\g (\pi)^{\leq r}]$ such that 
$$\forall v \in V, \quad F_X^{\pi, \leq r} (\widehat{v}_\pi)=(X \cdot \widehat{v}_\pi)^{\leq r}.$$
For all $\lambda_\pi \in \Co_\h^\g (\pi)$, one has 
$\delta_X^\pi (\lambda_\pi) =  \lim\limits_{\mathcal F}  \delta_X^{\pi, \leq r}(\lambda_\pi^{\leq r})$ and $\delta_X$ is continuous for the topology defined by the filtration ${\mathcal F}$.

\begin{notation}\label{definition delta_u}
If $u\in U(\g)$, its action on $\Co^\g_{ \h}{(\pi)}$ defines a differential operator on 
$\Co^\g_{ \h}{(\pi)}$ denoted $\delta_u^{\pi}$. If $r \in \N^*$, we define the differential operator 
$\delta_u^{\pi, \leq r}$ on  $\Co^\g_{ \h}{(\pi)}^{\leq r}$ as follows: 
$$\forall \lambda_\pi \in \Co^\g_{ \h}{(\pi)}^{\leq r}, \quad  
\delta_u^{\pi, \leq r}(\lambda_\pi):= \delta_u^\pi(\lambda_\pi )^{\leq r}.$$
$\delta_u^\pi $ is continuous with respect to the topology defined by the filtration $\F$. One has 
$$\delta _u^\pi  (\lambda_\pi)= \lim\limits_{\mathcal F}\delta_{u}^{\pi, \leq r} (\lambda_\pi^{\leq r})=
 \lim\limits_{\mathcal F} \left [ \delta^{\pi}_u (\lambda_\pi^{\leq r})\right ]^{\leq r}.$$

\end{notation}

\section{Restricted Lie superalgebra} 
Let $k$ be a field of characteristic $p>2$. A restricted Lie algebra $\g$ is a Lie algebra endowed with a  $p$-operation $(-)^{[p]} :\g \to \g$ , $X \mapsto X^{[p]}$ satisfying some special conditions (\cite{J}). A morphism of restricted Lie algebras is a map of Lie algebras preserving the $p$-operation.

\begin{definition}(\cite{P})
Let $k$ be a field of positive characteristic $p>2$.
 A Lie superalgebra $\g=\g_{\overline{0}}\oplus \g_{\overline{1}}$ is called restricted if 

\begin{itemize}
\item the Lie algebra $\g_{\overline{0}}$ is restricted;
\item the action of $\g_{\overline{0}}$ on $\g_{\overline{1}}$ defines a restricted morphism from 
$\g_{\overline{0}}$ to $\g l(\g_{\overline{1}})$.
\end{itemize}
A linear map   $f :\g \to \g^\prime$ is a morphism of restricted Lie superalgebras if 
$$\forall X \in \g_{\overline{0}}, \quad  f(X^{[p]}) =  f(X)^{[p]}.$$

\end{definition}

\begin{definition} Let $\mathfrak g$ be a restricted Lie superalgebra. A Lie subsuperalgebra $\mathfrak h$ of 
$\mathfrak g$ is a restricted subsuperalgebra of $\mathfrak g$ if 

(i) $\mathfrak h$ is a restricted Lie superalgebra.

(ii) The inclusion  map ${\mathfrak h} \to {\mathfrak g}$ is a morphism of restricted Lie superalgebras. 

\end{definition}

The restricted enveloping superalgebra of a restricted Lie superalgebra is defined as follows:
$$U^\prime({\mathfrak g})=\dfrac{U({\mathfrak g})}{\left ( X^p-X^{[p]}, \quad X \in \g_{\overline{0}}\right )}.$$

\begin{remark} If $\g$ is abelian and the $p$-operation is trivial, $U^\prime (\g)$ will be denoted $S^\prime (\g)$. Thus 
$$S^\prime({\mathfrak g})=\dfrac{S({\mathfrak g})}{\left ( X^p, \quad X \in \g_{\overline{0}}\right )}.$$
\end{remark}

Let $\g$ be a restricted $k$-Lie superalgebra and $\h$ a restricted Lie subsuperalgebra of $\g$. 
The Lie superalgebra $\h$ acts on $\g/\h$ by $ad_{\g/\h}$. Thus, the superalgebra $U(\h)$ acts on 
$\g/\h$ and, for $H\in \h_{\overline{0}}$, 
$$ad_{\g/\h}(H^p)=ad_{\g/\h}(H)^p.$$
As  $(ad_\g(H))^p=ad_\g(H^{[p]})$, one has 
$(ad_{\g/\h}(H))^p=ad_{\g/\h}(H^{[p]})$. As a consequence $ad_{\g/\h}$ is a representation of $U^\prime (\h)$ over 
$\g/\h$. 

The character $strad_{\g/\h}$ is well defined as a character of  $U^\prime(\h)$. \\

 We will now concentrate on the restricted crossed product defined by the restricted coinduced representation. 
 
 \begin{definition}

If $V$ is a $U^\prime ({\mathfrak h})$-module. We define its coinduced representation as 
$U^\prime ({\mathfrak g})$ acting on $Hom_{U^\prime({\mathfrak h})}(U^\prime ({\mathfrak g}),V)$ 
by the transpose of the right multiplication. 

\end{definition}

 For this section, we set 
$$\begin{array}{l}{\mathcal A}= Hom_{U^\prime({\mathfrak h})}(U^\prime ({\mathfrak g}),k), \\
\Co_{U^\prime (\mathfrak h)}^{U^\prime (\mathfrak g )}(V)=
Hom_{U^\prime({\mathfrak h})}(U^\prime ({\mathfrak g}),V).
\end{array}$$

${\mathcal A}$ is a local $k$-superalgebra with maximal ideal 
$${\mathcal M}=\{\eta \in {\mathcal A}, <\eta, 1>=0\}.$$

Moreover,  $\Co_{U^\prime (\mathfrak h)}^{U^\prime (\mathfrak g )}(V)$ is an
${\mathcal A}$-module and a $U^\prime (\g)$-module. 
It is a restricted crossed product (\cite{C4} for example). \\

The Poincaré-Birkhow-Witt theorem holds for restricted Lie  superalgebras (\cite{P} for example).

\begin{theorem}

Let $\g$ be a restricted Lie superalgebra. 
Suppose that $(e_i)_{i \in I}$ is an ordered  basis of ${\mathfrak g}_{\overline{0}}$ and $(\epsilon_j)_{j\in J}$ is an ordered basis of $\g _{\overline{1}}$. 

The monomials $e_{i_1}^{a_{i_1}}\dots e_{i_n}^{a_{i_n}}\epsilon_{j_1}...\epsilon_{j_t}$ with 
$i_1<\dots <i_n$, $ j_1<\dots <j_t$, $a_{i_j} \in [0,p-1]$  form a basis of the restricted enveloping algebra $U^\prime (\g)$.
\end{theorem}

\begin{notation}\label{restricted notation}
Let $\mathfrak h$ be a subrestricted Lie superalgebra of $\mathfrak g$. Given a supplement  ${\mathfrak p}$ of $\mathfrak h$ in $\mathfrak g$ and a basis $(e_1, \dots e_n, \epsilon_1, \dots , \epsilon_m)$ of it.
Set $(\underline{a}, \underline{\alpha})= (a_1, \dots ,a_n, \alpha_1, \dots , \alpha_m)\in [0,p-1]^{n,m}$. As before (\cite{C1}), the  map
 ${\mathcal J}_0 :  \Co_{U^\prime (\mathfrak h )}^{U^\prime (\mathfrak g)}(k) \to Hom_k \left ({\mathcal S}^\prime({\mathfrak p}), k \right )$
defined by: 
 $$\forall f \in \Co_{U^\prime (\mathfrak h )}^{U^\prime (\mathfrak g )}(k), \quad <e^{\underline{a}}\epsilon^{\underline{\alpha}}, {\mathcal J}_0(f)>= 
 <e^{\underline{a}}\epsilon^{\underline{\alpha}}, f>.$$
  is an isomorphism of superalgebras.

   If $(\pi, V)$ is a $U^\prime(\h)$-module, 
   we also define  an isomorphism 
 $J_\pi : \Co_{U^\prime({\mathfrak h})}^{U^\prime ({\mathfrak g})}(\pi) \to Hom \left (S^\prime({\mathfrak p} \right ), V)$ by: For all 
 $(\underline{a}, \underline{\alpha}) \in [0,p-1]^{n,m}$ and 
 $\eta \in \Co_{U^\prime ({\mathfrak h})}^{U^\prime({\mathfrak g})}(\pi)$, 
 $$<e^{\underline{a}}\epsilon^{\underline{\alpha}}, J_\pi (\eta)>= 
 <e^{\underline{a}}\epsilon^{\underline{\alpha}},\eta >.$$
 One has for any $f \in {\mathcal A}$ and $\eta \in \Co_{U^\prime ({\mathfrak h})}^{U^\prime({\mathfrak g})}(\pi)$,
 $J_\pi (f\eta)=J_0(f)J_\pi (\eta ).$

Let us introduce the elements $(\eta_1, \dots , \eta_n, \zeta_1, \dots , \zeta_m)$  of ${\mathcal A}$ defined by $\forall (a_1, \dots ,a_n, \alpha_1, \dots , \alpha_m)\in [0,p-1]^{n,m}$, 
$$\begin{array}{l}
< e_1^{a_1} \dots e_n^{a_n}\epsilon_1^{\alpha_1}\dots \epsilon_m^{\alpha_m}, \eta_i>=\delta_{0,a_1}\dots 
\delta_{0,a_{i-1}}\delta_{1,a_i}\delta_{0,a_{i+1}}\dots\delta_{a_n,0}\delta_{0,\alpha_{1}}\dots \delta_{0,\alpha_m};\\
< e_1^{a_1} \dots e_n^{a_n}\epsilon_1^{\alpha_1}\dots \epsilon_m^{\alpha_m}, \zeta_s>=
\delta_{0,a_{1}}\dots\delta_{a_n,0}\delta_{0,\alpha_{1}}\dots \delta_{0,\alpha_{s-1}}\delta_{1,\alpha_{s}}
\delta_{0,\alpha_{s+1}}\dots  \delta_{0,\alpha_m};\\
\end{array}$$
${\mathcal M}$ is the ideal of ${\mathcal A}$ generated by the elements $\eta_1, \dots , \eta_n, \zeta_1, \dots , \zeta_m$. 
 We set $$\Lambda=\eta_1^{p-1}\dots \eta_n^{p-1} \zeta_1 \dots \zeta_m.$$
If $(V,\pi)$ is a $U^\prime(\h)$-module and $v\in V$, we denote by $\hat{v}_\pi$ or  $\hat{v}$ the element of $Hom_{U^\prime (\h)}(U^\prime (\g), V)$ determined by 
$\forall (a_1, \dots ,a_n, \alpha_1, \dots , \alpha_m)\in [0,p-1]^{n,m}$, 
$$< e_1^{a_1} \dots e_n^{a_n}\epsilon_1^{\alpha_1}\dots \epsilon_m^{\alpha_m}, \widehat{v}_\pi>=
\delta_{0,a_1}\dots \delta_{0,a_n}\delta_{0,\alpha_{1}}\dots \delta_{0,\alpha_m}v.$$

If $u \in U(\g)$, we will denote by $\delta^\pi_u$ the coinduced action of $u$ on 
$\Co_{U^\prime (\h)}^{U^\prime ({\mathfrak \g})}(\pi)$. As before, we will write $\delta_u$ for $\delta^k_u$.

\end{notation}

 The following result was obtained in \cite{F-S} for $\g$  finite dimensional and totally even. \\

\begin{theorem}\label{Ind and  Coind}(\cite{F-S})
\begin{enumerate}
\item Let us set 
$$\Gamma^1_{\mathcal M}  \Co^{U^\prime (\g)}_{U^\prime (\h)}(\pi )=
\{\eta \in \Co^{U^\prime (\g)}_{U^\prime (\h)}(\pi ), \; \forall a \in {\mathcal M}, \; a\eta=0\}.$$

$\Gamma^1_{\mathcal M}{\mathcal A}$ is a (1,0) or (0,1) dimensional $k$-vector space with basis \linebreak 
$\Lambda=\eta_1^{p-1}\dots \eta_n^{p-1} \zeta_1 \dots \zeta_m$. The coinduced action  endows $\Gamma_{\mathcal{M}}{\mathcal A}$  with a $U^\prime (\h)$-module structure determined by the character $strad_{\g/\h}$.

\item Let $(\pi, V)$ be a restricted representation of $\h$. The $U^\prime (\h)$-module  $\Gamma^1_{\mathcal M} \Co^{U^\prime (\g)}_{U^\prime (\h)}(\pi )$ is isomorphic to 
$\Gamma^1_{\mathcal M} \Co^{U^\prime (\g)}_{U^\prime (\h)}(k )\otimes \pi.$

The map $$\begin{array}{rcl}
\Phi : \Ind_{U^\prime (\h)}^{U^\prime (\g)}\left (\Gamma^1_{\mathcal M}\Co^{U^\prime (\g)}_{U^\prime (\h)}(\pi ) \right ) 
& \to & \Co^{U^\prime (\g)}_{U^\prime (\h)}(\pi )\\
u\otimes_{U^\prime(\h)}(\Lambda  \otimes v) & \mapsto & \delta^\pi_{u} (\Lambda \hat{v}_\pi) 
\end{array}$$ is a $U^\prime (\g)$-isomorphism.
\end{enumerate}
\end{theorem}

{\it Proof of Theorem \ref{Ind and  Coind}:} 
\begin{enumerate}
\item The restriction of the coinduced action to $U^\prime (\h)$ preserves 
$\Gamma^1_{\mathcal M} \Co^{U^\prime (\g)}_{U^\prime (\h)}(k )$.
It is easy to see that  the superspace $\Gamma^1_{\mathcal M} \Co^{U^\prime (\g)}_{U^\prime (\h)}(k )$ has 
$\eta_1^{p-1}\dots \eta_n^{p-1} \zeta_1 \dots \zeta_m=\Lambda$ as a basis. 
Let us now check the following relation 
$$\forall H\in \h, \quad \delta_H(\Lambda )=strad_{\g/\h}(H) \Lambda.$$

If $H\in \h$, let us set 
$$[H,e_i]= \sum_{k=1}^n ad(H)_{k,i}e_k +  \sum_{s=1}^m ad(H)^{s,i}\epsilon_s \quad mod\; \h .$$ 
If $H \in \h_{\overline{1}}$, then $ad(H)_{i,i}=0$. 

Let $e^{\underline{a}}\epsilon^{\underline{\alpha}} \in U^\prime(\g)$ with $(\underline{a}, \underline{\alpha})\in [0,p-1]^{n,m}$. Let us write 
$$ e^{\underline{a}}\epsilon^{\underline{\alpha}}H=\sum_{(\underline{b},\underline{\beta}) \in [0,p-1]^{n,m}}f^{\underline{a}, \underline{\alpha}}_{\underline{b}, \underline{\beta}}(H)
e^{\underline{b}}\epsilon^{ \underline{\beta}}.$$

Let us denote by $\underline{p-1}$ (respectively $\underline{1}$) the element of $[0,p-1]$ 
(respectively $\{0,1\}^m$) whose components are all equal to $p-1$ (respectively $1$).
If $H \in \h$, the following equality holds: 
$$<f^{\underline{a},\underline{\alpha}}_{\underline{p-1},\underline{1}}(H)e^{\underline{p-1}}\epsilon^{\underline{1}}, \Lambda> = 
<e^{\underline{a}}\epsilon^{\underline{\alpha}}, \delta_H(\Lambda)>.$$
The coefficient  $f^{\underline{a},\underline{\alpha}}_{\underline{p-1},\underline{1}}(H)$ 
is zero if 
$e^{\underline{a}}\epsilon^{\underline{\alpha}}\neq e^{\underline{p-1}}\epsilon^{\underline{1}}$. Moreover, 
$$\begin{array}{rcl}
<e^{\underline{p-1}}\epsilon^{\underline{1}}, \delta_H(\Lambda )>&=&
f^{\underline{p-1},\underline{1}}_{\underline{p-1},\underline{1}}(H)<e^{\underline{p-1}}\epsilon^{\underline{1}}, \Lambda>\\
&=&-\left [ \sum_{i=1}^n (p-1)ad(H)_{i,i}- 
\sum_{s=1}^m ad(H)^{s,s}\right ]<e^{\underline{p-1}}\epsilon^{\underline{1}}, \Lambda>\\
&=&\left ( \sum_{i=1}^n ad(H)_{i,i}- \sum_{s=1}^m ad(H)^{s,s}\right )<e^{\underline{p-1}}\epsilon^{\underline{1}}, \Lambda>.\\
\end{array}$$

Thus 
$$\delta_H(\Lambda )=
strad_{\g/\h}(H)\Lambda.$$

\item 
One sees  easily  that the  $U^\prime (\h)$-module  
$\Gamma^1_{\mathcal M} \Co^{U^\prime (\g)}_{U^\prime (\h)}(\pi )$ is isomorphic to 
$\Gamma^1_{\mathcal M} \Co^{U^\prime (\g)}_{U^\prime (\h)}(k )\otimes \pi.$
Let us now show that the map $$\begin{array}{rcl}
\Phi : Ind_{U^\prime (\h)}^{U^\prime (\g)}\left [ \Gamma^1_{\mathcal M} \Co^{U^\prime (\g)}_{U^\prime (\h)}(\pi )\right ] 
& \to & \Co^{U^\prime (\g)}_{U^\prime (\h)}(\pi )\\
u\otimes_{U^\prime(\h)}(\Lambda  \otimes v) & \mapsto & \delta^\pi_{u}(\Lambda  \hat{v}_\pi) 
\end{array}$$
is well defined.
Let $v\in V$, $H \in \h$, $u,w \in U(\g)^\prime$.  On one hand, one has:
$$\begin{array}{rcl}
<w, \Phi_{uH \otimes (\Lambda \otimes v)}>&=& <wuH, \Lambda  \hat{v}_\pi >\\
&=&
<wu, \delta_H^{\pi}(\Lambda  \hat{v}_\pi)>\\
&=& <wu, \delta_H(\Lambda ) \hat{v}_\pi+ (-1)^{m\mid H\mid }\Lambda \delta_H^{\pi }(\hat{v}_\pi)>\\
&=& <wu, \strad_{\g/\h}(H)\Lambda  \hat{v}_\pi +(-1)^{m\mid H\mid } \Lambda \delta_H^{\pi }(\hat{v}_\pi)>\\
\end{array}$$ 
On the other hand:
$$\begin{array}{rcl}
<w, \Phi_{u \otimes H \cdot (\Lambda \otimes v)}>
&=& <wu,\strad_{\g/\h}(H)\Lambda  \hat{v}_\pi+ (-1)^{m\mid H\mid }\Lambda(\widehat{\pi(H)v})_\pi>\\

\end{array}$$ 

To finish the proof of Theorem \ref{Ind and  Coind}, we need to prove that 
$$<wu, \Lambda \delta_H (\hat{v}_\pi)>= <wu,\Lambda (\widehat{\pi(H)v})_\pi>.$$
But  one has 
$<wu, \Lambda >=0$ except if 
$wu= e_{1}^{p-1}\dots e_n^{p-1}\epsilon_1 \dots \epsilon_m.$ Thus, 
$$\begin{array}{rcl}
<wu, \Lambda \delta_H^{\pi}(\hat{v}_\pi)>&=&
<wu, \Lambda ><1,\delta_H (\hat{v}_\pi)>\\
&=&<wu, \Lambda ><H, \hat{v}_\pi>\\
&=&<wu, \Lambda >\pi(H)(<1, \hat{v}_\pi>)\\
&=&<wu, \Lambda>\pi(H)(v)\\
&=& <wu,\Lambda><1,(\widehat{\pi(H)v})_\pi>\\
 &=&<wu,\Lambda(\widehat{\pi(H)v})_\pi>.
 \end{array}$$

 Set $\check{\underline{a}}=\underline{p-1}-\underline{a}$ and $\check{\underline{\alpha}}=\underline{1}-\underline{\alpha}.$
For any $\eta \in  \Co^{U^\prime(\g)}_{U^\prime(\h)} (\pi)$, one has 
$$\eta= \Phi \left (\sum_{\underline{a}, \underline{\alpha}}
\dfrac{1}{<e^{\check{\underline{a}}}\e^{\check{\underline{\alpha}}}e^{\underline{a}}e^{\underline{\alpha}},\Lambda>}
e^{\underline{a}} \e^{\alpha} \otimes
\left [\Lambda \otimes \widehat{<e^{\check{\underline{a}}}\e^{\check{\underline{\alpha}}},\eta>_\pi}\right ]\right ).\Box 
$$
\end{enumerate}

\medskip

We will now construct a $U^\prime (\g)$-invariant duality between 
${\Co}_{ U^\prime (\h)}^{U^\prime (\g)}(\pi)$ and  
${\Co}_{ U^\prime (\h ) }^{U^\prime (\g)} (\pi^* \otimes  \Pi^n Ber ( (\g /\h)^*)$ This will rely on an integration argument on ${\Co}_{ U^\prime (\h)}^{U^\prime (\g)}(\Pi^n Ber (\g/\h)^*)$. For that purpose, we will identify 
${\Co}_{ U^\prime (\h)}^{U^\prime (\g)}(\Pi^n Ber (\g/\h)^*)$ with  the volume forms on ${\mathcal A}$.

\medskip

Let us first describe  the action of an element $X \in {\mathfrak g}$ on  a coinduced module and introduce some notation.

 \begin{notation}\label{notation delta_X}  
  The derivations $(\partial_1, \dots , \partial_n, \overline {\partial}_1, \dots , \overline{\partial}_s)$ defined by 
 $$\partial_i(\eta_j)=\delta_{i,j},\quad \partial_i(\zeta_t)=0, \quad  \overline{\partial}_s(\eta_j)=0, 
 \quad \overline{\partial}_s(\zeta_t)=\delta_{s,t}$$ 
 form a basis of the ${\mathcal A}$-module 
 $Der ({\mathcal A})$.

 Let $X \in {\mathfrak g}$. Denote by $\delta_X$ the transpose of the right multiplication on 
 $U^\prime({\mathfrak g})$. It is easy to see that, if $k\geq 1$, then $\delta_X$ sends   ${\mathcal M}^{k}$ to ${\mathcal M}^{k-1}$. 
 It defines a derivation of ${\mathcal A}$ that will be written 
$$\delta_X=\sum_{i=1}^n f_{i}(X)\partial_{i}+ \sum_{k=1}^m g_s(X) \overline{\partial}_s.$$
We can be more precise. 
\begin{equation}\label{explicit delta_X}
\begin{array}{rcl}
\delta_X&=&\sum_{i=1}^n <X, \eta_i>\partial_i + \sum_{s=1}^m <X, \zeta_s>\overline{\partial}_s \\
&+&\sum_{i=1}^n \sum_{\underline{a}, \underline{\alpha}}<e^{\underline{a}}\epsilon^{\underline{\alpha}}X, \eta_i>\partial_i 
+\sum_{s=1}^m \sum_{\underline{a}, \underline{\alpha}}<e^{\underline{a}}\epsilon^{\underline{\alpha}}X, \zeta_s>\overline{\partial_s}.
\end{array}
\end{equation}

The element $X \in \g$ defines a differential operator of degree 1 ,  $\delta_X^\pi$ (denoted $\delta_X$ if there is no ambiguity), of $\Co^{U^\prime (\g)}_{ U^\prime (\h)}{(\pi)}$. If 
 $F^\pi _X$ denotes  the element of $End_A[\Co_{U^\prime (\h)}^{U^\prime (\g)} (\pi)]$ defined by 
$$\forall v \in V, \quad F_X^\pi (v)=X\cdot v.$$
Then 
\begin{equation}\label{explicit delta_X^pi}
\delta_X^\pi=F^\pi _X + \sum_{i=1}^n  f_{i}(X)\partial_{i}^\pi
+\sum_{s=1}^m g_s (X)\overline{\partial_s}^\pi.\end{equation}

\end{notation}

\begin{proposition}\label {Omega and coinduced}(\cite{C1})
Set $Der({\mathcal A})^*= Hom_{\mathcal A}\left ( Der({\mathcal A}), {\mathcal A} \right )$ and 
  $\Omega:= Ber\left (Der ({\mathcal A})^* \right )$.
\begin{enumerate}
\item  The ${\mathcal A}$-module $\Omega$ is endowed  with  a 
$U^\prime(\g)$-module by the operations:
$$\begin{array}{l}
  \forall X \in \g,  \quad \forall \omega \in \Omega, \quad 
X\cdot \omega =L_{\delta_X}(\omega).
\end{array}$$

\item 
The ${\mathcal A}$-module $\Omega $ is a free  ${\mathcal A}$-module of rank one with basis \linebreak 
$\omega_{\underline{e}}=\Pi \partial_{1}^*\dots \Pi \partial_n^*\overline{\partial_{1}}\dots \overline{\partial_m}$.

The map $\sigma : \g \to Der({\mathcal A}), \quad X \mapsto \delta_X$ is a morphism of $\g$-modules. 
It induces an isomorphism of $\h$-modules from $Ber \left ( \g/\h \right )^*$ to
$ \dfrac{\Omega}{{\mathcal M}\Omega}$
that sends $\Pi e_1^*\dots \Pi e_n^*\epsilon_1 \dots \epsilon_m$ to 
$\overline{\omega_{\underline{e}}}:=\Pi \partial_1^*\dots \Pi \partial_n^*\overline{\partial_1}\dots \overline{\partial _m}$ mod ${\mathcal M}$. 
The $k$-vector superspace $Ber \left [ \dfrac{\Omega}{{\mathcal M}\Omega}\right ]$  
 is naturally endowed with a $U^\prime (\h)$-module structure which is given by the character $-strad_{\g/\h}$.

\item The map 
\begin{equation}\label{chi_r}
\begin{array}{rcl} 
\chi : \Omega  & \to & \Co_{U^\prime (\h)}^{U^\prime(\g)} (\dfrac{\Omega}{{\mathcal M}\Omega})
\\
\omega & \mapsto & [X_{1}\dots X_{t} \in U^\prime (\g)\mapsto 
(\delta_{X_1}\dots \delta_{X_t} \cdot \omega)\; mod \;{\mathcal M}]
\end{array}
\end{equation}
is an isomorphism of ${\mathcal A}-U^\prime (\g)$-modules from $\Ber(Der \left ( {\mathcal A})^*\right )$  to 
${\Co}_{U^\prime (\h)}^{U^\prime (\g)} \left (\Ber (\g/\h )^*\right )\simeq 
{\Co}_{U^\prime (\h)}^{U^\prime (\g)} (\Pi^{n+m}k_{-strad_{\g/\h}})$.

Moreover, there exists an invertible element $g$ 
 of ${\mathcal A}$ such that (using Notation (\ref{restricted notation}))
 $\chi \left [ \omega_{\underline{e}} \right ]=g \widehat{\overline{\omega_{\underline{e}}}}.$\\
\end{enumerate}
\end{proposition}

{ \it Proof of Proposition \ref{Omega and coinduced}:} 
\begin{enumerate}
\item The first assertion is an easy computation.

\item As for any $H \in \h$, $\sigma(H) ({\mathcal M})\subset {\mathcal M}$,  the morphism $\sigma$ induces an isomorphism of $\h$-modules 
$$\begin{array}{rcl}
\g/\h & \to & \dfrac{Der({\mathcal A})}{{\mathcal M}Der({\mathcal A})}
\end{array}$$
that sends $e_i$ to $\delta_{e_i}=\partial_i$ mod ${\mathcal M}$ (using \ref{explicit delta_X}). 
It induces an isomorphism of $\h$-modules from $Ber \left ( \g/\h \right )^*$ to
$ \dfrac{\Omega}{{\mathcal M}\Omega}$.

\item It is easy to check that $\chi$ is a morphism of $U^\prime ( \g )$-modules and of ${\mathcal A}$-modules.  
Let us now show that it is an isomorphism. 
It is  an isomorphism modulo ${\mathcal  M}$ as  
$\chi (\omega_{\underline{e}})(1)=\omega_{ \underline{e}} $ mod ${\mathcal M}$. 

Moreover, $\Omega$ and 
 $\Co_{U^\prime(\h)}^{U^\prime ( \g)}(\Omega/{\mathcal M}\Omega)$ 
 are free ${\mathcal A}$-modules of dimension 1 and $\chi $ sends a basis of $\Omega$ to a basis of 
 $\Co_{U^\prime (\h)}^{U^\prime (\g)}(\Omega/{\mathcal M}\Omega)$. 
 Thus, there exists an invertible element $g$ 
 of ${\mathcal A}$ such that 
 $$\chi \left [ \omega_{\underline{e}} \right ]=g\widehat{\overline{\omega_{\underline{e}}}}.\Box$$

\end{enumerate}

\bigskip

The choice of basis of a supplement of $\h$ in $\g$ defines coordinates on ${\mathcal A}$ (as in Notation \ref{restricted notation}) and a basis 
$\omega_{\underline{e}}=\Pi \partial_1^*\dots \Pi \partial_n^*\overline{\partial_1}\dots \overline{\partial _m}$
 of $\Omega$.
Using Proposition \ref {Omega and coinduced}, we define a map 
$$\Psi : {\Co}_{ U^\prime (\h)}^{U^\prime (\g)}(\pi) \otimes 
{\Co}_{ U^\prime (\h ) }^{U^\prime (\g)} (\pi^* \otimes  \Pi^n Ber ( (\g /\h)^*)\simeq 
{\Co}_{ U^\prime (\h)}^{U^\prime (\g)}(\pi) \otimes\left ({\Co}_{ U^\prime (\h ) }^{U^\prime (\g)} (\pi^*) \otimes_{\mathcal A}  \Pi^n\Omega\right ) \to k$$
$$\Psi \left [ \sum_{(\underline{a}, \underline{\alpha})} \eta^{\underline{a}} \zeta^{\underline{\alpha}}
\widehat {v}_{\underline{a}, \underline{\alpha}}, 
 \sum_{(\underline{b}, \underline{\beta})} \eta^{\underline{b}}\zeta^{\underline{\beta}}
 \widehat{v}^*_{\underline{b}, \underline{\beta}}\otimes \Pi^n \omega_{\underline{e}}\right ]=
 \sum_{a_{i}+b_{i}=p-1, \; \alpha_s+\beta_s=1} (-1)^{Inv(\underline{\alpha}, \underline{\beta})}
 (-1)^{\mid v_{\underline{a}, \underline{\alpha}}\mid \mid \beta \mid }
 <v_{\underline{a}, \underline{\alpha}}, v^*_{\underline{b}, \underline{\beta}}>
$$
where $\underline{a}, \underline{b}\in [0,p-1]^n$, $\underline{\alpha}, \underline{\beta }\in \{0,1\}^m$ and $\zeta^{\underline{\alpha}}\zeta^{\underline{\beta}}=(-1)^{Inv(\underline{\alpha}, \underline{\beta})}\zeta_1 \dots \zeta_m$.

If $\lambda \in \Co_{U^\prime (\h)}^{U^\prime (\g)}(\pi)$ and $\lambda^* \in \Co_{U^\prime (\h)}^{U^\prime (\g)}(\pi^*)$, $\Psi$ can be written as follows:
$$\Psi \left (\lambda, \lambda^* \otimes  {\Pi^ n} \omega_{\underline{e}}\right )=
(-1)^{\frac{m(m-1)}{2}}\dfrac{1}{(p-1)!^n}
<e_1^{p-1}\dots e_n^{p-1}\epsilon_1\dots \epsilon_m, <\lambda,\lambda^*>>
$$
where $<\lambda,\lambda^*>$ is the element of ${\mathcal A}$ defined by 
$$\forall u \in U(\g), \quad <u, <\lambda, \lambda^*>>=\sum < <u_{(1)},\lambda>, <u_{(2)},\lambda^*>>(-1)^{\mid u_{(2)}\mid \mid \lambda \mid}.$$

 \begin{corollary}\label{definition of Psi}
 The map $\Psi$ defines a non degenerate $U^\prime({\mathfrak g})$-invariant duality between 
 ${\Co}_{U^\prime ( \h)}^{U^\prime (\g)}(\pi)$ and  
$ {\Co}_{U^\prime ( \h) }^{U^\prime (\g)} (\pi^* \otimes \Pi^n Ber [(\g /\h)^*] )$. 

\end{corollary}

{\it Proof of Corollary \ref{definition of Psi}:}

We keep the same notation as above. 
We are in the situation where $Der({\mathcal A})$ is a finite dimensional free 
${\mathcal A}$-module of dimension $(n,m)$ with basis 
$( \partial _1 , \dots ,  \partial_n, \overline{\partial_1}, \dots ,\overline{\partial_m})$. 
The Lie derivative is defined on $\Ber(Der \left ( {\mathcal A})^*\right )$.

Let $\lambda \in \Co_{U^\prime(\h)}^{U^\prime ({\g})}(\pi)$, 
$\lambda^* \in \Co_{U^\prime(\h)}^{U^\prime ({\g})}(\pi^*)$, we need to show that 
$$\Psi \left [ L_{\delta_X}\left (<\lambda, \lambda^* >\Pi^n\omega_{\underline{e}}\right ) \right ]=0.$$
This relation follows from a computation using the explicit formula for $\delta_X$ (\ref{explicit delta_X}).

The maps $\Phi$ and $\Theta$ being isomorphisms of $U^\prime (\g)$-modules,  it is also an easy consequence of the relation between $^t\Phi$,  $\Theta$ and $\Psi$ proved in next Proposition \ref{comparison of Phi and Psi}. That is why, we won't reproduce the computation.$\Box $\\

It is well known that the map 
\begin{equation}\label{definition Theta}
\begin{array}{rcl}
\Theta : Coind_{U^\prime (\h)}^{U^\prime (\g)} (\pi^*)&\to &
 {\mathcal I}nd_{U^\prime (\h)}^{U^\prime (\g)} (\pi)^*\\
\lambda & \mapsto & \left [ u \otimes_{U^\prime (\h)} v\mapsto <\lambda(\check{u}),v>\right ]\\
\end{array}
\end{equation}
 is an isomorphism of $U^\prime (\g)$-modules. 
We will now make explicit the map 
$$ ^t{\Phi}^{-1}\circ \Theta : \Co_{U^\prime (\h)}^{U^\prime(\g)}(\pi ^*\otimes \Gamma_{\mathcal M}^1{\mathcal A}) \to \left [ \Co_{U^\prime (\h)}^{U^\prime(\g)}(\pi)\right ]^*.$$

\begin{proposition}\label{comparison of Phi and Psi}

\begin{enumerate}
\item The map 
\begin{equation}\label{Identification}
\begin{array}{rcl}
\Pi^n Ber(\g/\h) & \to &  \Gamma^1_{\mathcal M}\mathcal A \\
e_1 \dots e_n \epsilon_1^*\dots \epsilon_m^* & \to & \eta_1^{p-1} \dots\eta_n^{p-1} \zeta_1 \dots \zeta_m 
\end{array}
\end{equation} 
is an isomorphism of $U^\prime (\h)$-modules.
\item 

Let $\Phi : Ind_{U^\prime(\h)}^{U^\prime(\g)}( \pi \otimes \Gamma^1_{\mathcal M}\mathcal A) \to \Co_{U^\prime(\h)}^{U^\prime (\g)}(\pi )$ 
the isomorphism constructed in  Theorem \ref{Ind and Coind} and 
$^{t}\Phi : \Co_{U^\prime(\h)}^{U^\prime(\g)}(\pi )^* \to 
Ind_{U^\prime(\h)}^{U^\prime(\g)}(  \pi \otimes \Gamma^1_{\mathcal M}\mathcal A)^*$ its transpose. 
Denote by $\Psi^\natural : \Co_{U^\prime (\h)}^{U^\prime(\g) } \left (\pi ^*\otimes \Pi^n Ber(\g /\h)^*\right ) \to \left [ \Co_{U^\prime (\h)}^{U^\prime(\g)}(\pi)\right ]^*$ 
the monomorphism determined by $\Psi
$. If we identify  $\Pi^n Ber(\g/\h)$ with $\Gamma^1_{\mathcal M}{\mathcal A}$ by the map (\ref{Identification}), the following equality   holds (see \ref{definition Theta} for the definition of the map $\Theta$)
$$^t \Phi  \circ  \Psi ^\natural =\Theta .$$ 
\end{enumerate}
\end{proposition}

{\it Proof of Proposition \ref{comparison of Phi and Psi}}:

\begin{enumerate}
\item $ \partial_1 \dots  \partial_n \overline{\partial_1}^*\dots \overline{\partial_m}^* $ is a basis of the ${\mathcal A}$-module 
$ \Pi^n Ber \left ( Der ({\mathcal A}) \right )$. From (\ref{explicit delta_X}), 
$ \partial_1 \dots  \partial_n \overline{\partial_1}^*\dots \overline{\partial_m}^* $ mod ${\mathcal M} $ can be identified to $e_1\dots e_n \epsilon_1^*\dots \epsilon _m^*$. 
The assertion follows from Theorem \ref{Omega and coinduced}.

\item 
Let $\omega_{\underline{e}}=\Pi \partial_1^* \dots \Pi \partial_n^* \overline{\partial_1}\dots \overline{\partial_m} \in Ber \left [Der ({\mathcal A})^* \right ]$, 
$\overline{\omega_{\underline{e}}}=\Pi e_1^*\dots \Pi e_n^*\epsilon_1 \dots \epsilon_m \in Ber \left [Der ({\mathcal A}) \right ]^*/{\mathcal M}Ber \left [Der ({\mathcal A}) \right ]^* \simeq Ber \left [ \g/\h\right ]^*$ and 
$\Pi^n\widehat{\overline{\omega_{\underline{e}}}}$ the basis of $\Co_{U^\prime(\h)}^{U^\prime (\g)}( \Pi^n Ber (\g/\h )^*)$ it defines (see Notation \ref{restricted notation} and Theorem \ref{Omega and coinduced}). 
Let $u \in U(\g)$, $v \in V$, $f\in {\mathcal A}$. 
On one hand, we have:
$$\begin{array}{rcl}
<u\otimes_{U^\prime (\h)} v \Lambda, \Theta (f\hat{v^*}\otimes   \Pi^n\widehat{\overline{\omega_{\overline e}}})>&=&  <\check{u},f> <v,v^*> <\Lambda, \Pi^n\overline{\omega_{\overline e}}>(-1)^{m(\mid f \mid +\mid v^*\mid ) + \mid f\mid \mid v \mid }\\
&=& <\check{u},f> <v,v^*>(-1)^{m(\mid f \mid +\mid v^*\mid ) + \mid f\mid \mid v \mid }.
\end{array}$$

On the other hand, we have:

$$\begin{array}{rcl}
<u\otimes v \Lambda,    {^t \Phi }\circ \Psi^\natural \left (f \hat{v^*} \otimes \Pi^n  \widehat{\overline{\omega_{\overline e}}} \right )> &=& 
\Psi \left [ \delta_{u} (  \widehat{v} \Lambda), f\widehat{v^*}\otimes \Pi^n\widehat{\overline{\omega_{\underline{e}}}} \right ]\\
&=& \Psi \left [ \widehat{v}, \Lambda \delta_{\check{u}}(f \widehat{v^*} \otimes \Pi^n \widehat{\overline{\omega_{\underline{e}}}}),   \right ]
(-1)^{\mid u \mid ( \mid v\mid +m)}\\
&=& \Psi \left [ \widehat{v}, \Lambda  \delta_{\check{u}}(f) \widehat{v^*}\otimes \Pi^n \widehat{\overline{\omega_{\underline{e}}}}   \right ]
(-1)^{\mid u \mid ( \mid v\mid +m)}\\
&=&\Psi \left [ \hat{v},\Lambda <\check{u},f>\widehat{v^*}\otimes \Pi^n \widehat{\overline{\omega_{\underline{e}}}}\right ]
(-1)^{\mid u \mid ( \mid v\mid +m  )}\\
&\underset{Th. \ref{Omega and coinduced}}{=}&\Psi \left [ \hat{v},\Lambda <\check{u},f>\widehat{v^*}<1,g>\Pi^n\omega_{\underline{e}}\right ]
(-1)^{\mid u \mid ( \mid v\mid +m  )}\\
&=&<\check{u}, f><v, v^*>(-1)^{\mid u \mid ( \mid v\mid +m  )+m \mid v^* \mid }\\
&=&<u\otimes_{U^\prime(\h)} ( v \Lambda ), 
\Theta (f  \widehat{v^*} \otimes \Pi^n \widehat{\overline{\omega_{\overline e}}})>
\end{array}$$
The last equality follows from the fact that if $<\check{u}, f>\neq 0$, then $\mid f \mid =\mid u \mid $. 
\end{enumerate}

\section{Applications} 
Let $k$ be a field of characteristic $p>2$. 
Let $\g$ be a Lie superalgebra and $\h$ a be subsuperalgebra of finite codimension. Let $\mathfrak p$ be a supplement of $\h$ in $\g$ so that $\g=\h\oplus {\mathfrak p}$. 
 Let 
$\underline{e}:= (e_1, \dots , e_n, \epsilon_1, \dots , \e_m )$ be a basis of $\mathfrak p$.  For all $i\in[1,n]$, assume given 
$$z_i=e_{i}^{p^{t_i}}-\sum_{j=0}^{{t_i}-1}x_{i,j}e^{p^j}$$
 with $x_{i,j }\in k$.

\begin{example}
 If $\h=\{0\}$ so that $\g$ is finite dimensional, one can construct $(z_1,\dots,z_n) $ in the center of $U(\g)$ (\cite{F-S1}).

 \end{example}


Let  $Z=\left ( e^{p^{t_i}}- \sum_{j=0}^{t_i-1} x_{i,j} e^{p^j}, \quad i \in [1,n]\right )$ be the two sided ideal of $U(\g)$ generated by the elements 
$\{z_i=e^{p^{t_i}}- \sum_{j=0}^{t_i-1} x_{i,j} e^{p^j}, \quad i \in [1,n]  \}$. As the  $z_i$'s are primitive elements of $U(\g)$, $Z$  is a coideal. 
We will be interested in the following quotient of $U(\g)$.
$$U_Z({\mathfrak g})=\dfrac{U({\mathfrak g})}
{\left ( e^{p^{t_i}}- \sum_{j=0}^{t_i-1} x_{i,j} e^{p^j}, \quad i \in [1,n], \; j \in \N\right )}.$$
As $\check{z_i}=-z_i$, the superbialgebra $U_Z(\g)$ is a cocommutative Hopf superalgebra.
Poincaré-Birkhoff Witt Theorem holds for $U_Z(\g)$.

\begin{theorem}

We keep the same notation as above.
\begin{enumerate}\item The monomials $z_{1}^{c_{1}}\dots z_{n}^{c_{n}}e_{1}^{a_{1}}\dots e_{n}^{a_{n}}\epsilon_{j_1}...\epsilon_{j_t}$ with 
$ j_1<\dots <j_t$, $a_i \in [0,p^{t_i}-1]$, $c_{k}\in \N$  form a basis of 
$U(\h)$-module $U(\g)$.
\item The monomials $e_{1}^{a_{1}}\dots e_{n}^{a_{n}}\epsilon_{j_1}...\epsilon_{j_t}$ with 
 $ j_1<\dots <j_t$, and $a_{i} \in [1,p^{t_i}-1]$  form a basis of 
$U(\h)$-module $U_Z(\g)$.
\end{enumerate}
\end{theorem}

 We will now concentrate on generalized  coinduced representations. 
 
 \begin{definition}

If $V$ is a $U({\mathfrak h})$-module. We define its generalized coinduced representation as 
$U_Z ({\mathfrak g})$ acting on ${\mathcal A}_Z=Hom_{U({\mathfrak h})}(U_Z ({\mathfrak g}),V)$ 
by the transpose of the right multiplication. 

\end{definition}

Denote by $P_Z$ the set of primitive elements of $U_Z(\g)$. As before, we set for $i \in [1,n]$ and $j \in [0,t_i-1]$, $e_{i,j}=e_i^{p^j}$. One has 
$$P_Z=Vect\{e_{i,j}, \epsilon_s\quad {\rm for} \; i \in [1,n]\;  \rm{and} \; j \in [0,t_i-1], s \in [1,m]\}.$$
Then $P_Z$ is a restricted Lie subsuperalgebra. 
As $U_Z(\g)$ is a  cocommutative Hopf superalgebra generated by  $P_Z$, from the  Milnor-Moore  theorem (\cite{M-M}), the Hopf algebras $U_Z(\g)$ and $U^\prime(P_Z)$ are isomorphic. We can apply the result of the previous section and we get the following corollary. As usual we denote by $(\eta_{i,j}, \zeta_s)$ the basis if $(P_Z)^*$ dual to 
$\left (e_{i,j}, \epsilon_s\quad {\rm for} \; i \in [1,n]\;  \rm{and} \; j \in [0,t_i-1], s \in [1,m]\right )$. 

\begin{corollary}\label{Ind and  Coind over A_Z}
\begin{enumerate}
\item Set 
$\displaystyle \Lambda= \left ( \prod_{i=1}^n \prod_{j=0}^{p^{t_i-1}}\eta_{i,j}^{p-1}\right )\zeta_1 \dots \zeta_m$. 
The one dimensional space $k\Lambda$  is endowed with the following $\h$-module structure  
$$\forall H\in \h, \quad \delta_H(\Lambda )=strad_{\g/\h}(H) \Lambda.$$
\item Let $(\pi, V)$ be a  representation of $\h$. 
The map $$\begin{array}{rcl}
\Phi : Ind_{U (\h)}^{U_Z (\g)}( \Pi^m k_{\strad_{\g/\h}} \otimes \pi)  & \to & \Co^{U_Z (\g)}_{U (\h)}(\pi )\\
u\otimes_{U(\h)}(\Lambda  \otimes v) & \mapsto & \delta^\pi_{u} (\Lambda \hat{v}_\pi) 
\end{array}$$ is a $U_Z (\g)$-isomorphism.
\item Set $\Omega:= Ber\left ( Der({\mathcal A}_Z)\right )$ and $\tau=t_1+\dots +t_n$. Using Proposition \ref {Omega and coinduced}, we define a map 
$$\Psi : {\Co}_{ U (\h)}^{U_Z(\g)}(\pi) \otimes 
{\Co}_{ U (\h ) }^{U_Z (\g)} (\pi^* \otimes  \Pi^\tau Ber  (P_Z /\h)^*)\simeq 
{\Co}_{ U (\h)}^{U_Z (\g)}(\pi) \otimes\left ({\Co}_{ U (\h ) }^{U_Z (\g)} (\pi^*) \otimes_{\mathcal A}  \Pi^\tau \Omega\right ) \to k.$$
If $\lambda \in \Co_{U (\h)}^{U_Z(\g)}(\pi)$ and $\lambda^* \in \Co_{U (\h)}^{U_Z (\g)}(\pi^*)$, $\Psi$ can be written as follows:
$$\Psi \left (\lambda, \lambda^*\omega_{\underline{e}}\right )=
(-1)^{\frac{m(m-1)}{2}}\dfrac{1}{  (p-1)!^{\tau}}
<e_{1}^{p^{t_1}-1}\dots e_{n}^{p^{t_n}-1}\epsilon_1\dots \epsilon_m, <\lambda,\lambda^*>>
$$
where $<\lambda,\lambda^*>$ is the element of ${\mathcal A}$ defined by 
$$\forall u \in U(\g), \quad <u, <\lambda, \lambda^*>>=\sum < <u_{(1)},\lambda>, <u_{(2)},\lambda^*>>(-1)^{\mid u_{(2)}\mid \mid \lambda \mid}.$$
 
 The map $\Psi$ defines a non degenerate $U_Z({\mathfrak g})$-invariant duality between 
 ${\Co}_{U ( \h)}^{U_Z (\g)}(\pi)$ and  
${\Co}_{U ( \h) }^{U_Z (\g)} (\pi^* \otimes \Pi^m k_{-strad_{\g/\h}} )$. 
\item Denote by $\Theta$ the following isomorphism of $U_Z(\g)$-modules
$$\begin{array}{rcl}
\Theta : Coind_{U (\h)}^{U_Z (\g)} (\pi^*)&\to &
 {\mathcal I}nd_{U (\h)}^{U_Z (\g)} (\pi)^*\\
\lambda & \mapsto & \left [ u \otimes_{U (\h)} v\mapsto <\lambda(\check{u}),v>\right ].
\end{array}$$ 
Denote by 
$\Psi^\natural : \Co_{U (\h)}^{U_Z(\g)}(\pi ^*\otimes \Pi^m k_{-strad_{\g/\h}}) \to 
\left [ \Co_{U (\h)}^{U_Z(\g)}(\pi)\right ]^*$ 
the isomorphism  determined by $\Psi$.
The following   equality holds 
$$^t \Phi  \circ  \Psi ^\natural =\Theta.$$ 

\end{enumerate}
\end{corollary}

We will now study  a duality property for  the kernel of  coinduced representations. For that purpose, for a fixed $r$ in  $\N$. We will use the previous construction in the case where 
$z_i=e_i^{p^r}$. Thus we construct the restricted enveloping superalgebra 
$$U(\g)^{\leq r}=\dfrac{U(\g)}{\left (e_1^{p^{r+1}}, \dots, e_n^{p^{r+1}}\right )}.$$
$U(\g)^{\leq r}$ is isomorphic to ${\mathcal F}_rU(\g)$ as $U(\h)$-modules (see Notation \ref{usual notation}) and coalgebras.
One has 
$$U(\g)=\displaystyle{\underset{}{\lim_{\rightarrow}}}U(\g) ^{\leq r}=
\displaystyle{\underset{}{\lim_{\rightarrow}}}{\mathcal F}_r U(\g).$$

 \begin{remark}\label{preliminary remark}
 If $u\in U(\g)^{\leq r_0}$, it defines a differential operator on $\Co_{U(\h)}^{U (\g)^{\leq r}}(\pi)$ that can be identified with   $\delta_{u}^{\pi , \leq r}$ (see Notation \ref{definition delta_u}) if $r\geq r_0$. 
 \end{remark}

 \begin{theorem} \label{annihilators duality}
 Let $k$ be a field of characteristic $p>2$. Let $\g$ be a Lie $k$-superalgebra and $\h$ a Lie subsuperalgebra of $\g$. Let $(\pi,V)$ be a representation of $\h$. Denote by $I_\pi \subset U({\mathfrak g})$ the kernel of the representation $\Co_\h^\g (\pi)$. Assume that the $k$-vector space 
 ${\g}/{\h}$ is finite dimensional with odd dimension $m$.  Then 
 $$I_{\pi}=\check{I}_{\pi^*\otimes \Pi^{m}k_{-strad_{\g /\h}}}=\check{I}_{\pi^*\otimes k_{-strad_{\g /\h}}}$$
 where $k_{-strad_{\g /\h}}$  is the one dimensional representation of $\h$ defined by the character $-strad_{\g /\h}$.
 \end{theorem}
\begin{remarks}
Theorem \ref{annihilators duality} is proved in \cite{Du} for Lie algebras in any characteristic but with the assumption $\g$ finite dimensional. 
It is proved in \cite {C1} in the setting of Lie superalgebras for $\g/\h$ finite dimensional but only in characteristic $0$.
\end{remarks}
{\it Proof of Theorem  \ref{annihilators duality}:} 

First, the equality  
$I_{\pi^*\otimes \Pi^{m}k_{-strad_{\g /\h}}}=I_{\pi^*\otimes k_{-strad_{\g /\h}}}$ holds  as  a consequence of  the following remarks:

\begin{itemize}
\item If $Ann(\pi) \subset U(\h)$ is the annihilator of the representation $(\pi,V)$,  the annihilator of the $U(\g)$-module 
$ \Co_\h^\g (\pi )$ is the largest two sided ideal contained in  $Ann(\pi)U(\g)$ (\cite{Di}).
\item If $\chi$ is a character of $U(\h)$, then $\chi$ vanishes  on any  element of 
$U(\h)_{\overline{1}}$. Then, it is easy to see that  $U(\h)$-modules $\Pi (V\otimes k_\chi)$ and 
$V\otimes \Pi k_{\chi}$ have the same annihilator.

\end{itemize}

Thus we will prove the equality $I_{\pi}=\check{I}_{\pi^*\otimes \Pi^{m}k_{-strad_{\g /\h}}}$.

To simplify notation in the proof, we set $\tilde{\pi}:= \pi^*\otimes \Pi^{m}k_{-strad_{\g /\h}}$
Let $u \in {\mathcal F}_{r_0}U(\g)$. It defines a differential operator in $\Co_\h^\g (\pi)$ denoted $\delta_{u}^\pi$. If $r\geq r_0$, it also defines a differential operator in 
$\Co_{U(\h)}^{U (\g)^{\leq r}}(\pi)$ that we identify with  $\delta_{u}^{\pi , \leq r}$ (see Remark \ref{preliminary remark}).
From Corollary \ref{Ind and  Coind over A_Z}, one has 
$$\delta_{u}^{\pi , \leq r}=0 \Longleftrightarrow \delta_{\check {u}}^{\tilde{\pi} , \leq r}=0.$$
Then 
$$\begin{array}{rcl}
u \in I_\pi & \Longleftrightarrow &  \delta_{u}^{\pi }=0\\
& \Longleftrightarrow & \forall r \in \N, \quad \delta_{u}^{\pi , \leq r}=0\\
&\Longleftrightarrow & \forall r \in \N, \quad\delta_{\check {u}}^{\tilde{\pi} , \leq r}=0\\
&\Longleftrightarrow & \delta_{\check {u}}^{{\tilde{\pi}}}=0.\Box
\end{array}$$

\section{Declarations}

{\bf Ethical approval: } 
Non applicable 

{\bf Competing interest:  } 
The author declares no conflict of interest.

{\bf Authors' contributions: }
Non applicable.

{\bf Funding: } 
Non applicable.

{\bf Availability of data and materials: } 
Non applicable.

\vspace{1em}
Sophie Chemla

Sorbonne Université, Université  Paris Cité, CNRS, IMJ-PRG, F-75005 Paris

Email address: sophie.chemla@sorbonne-universite.fr

\end{document}